\documentclass[a4paper,12pt]{amsart}
\usepackage[utf8]{inputenc}
\usepackage{amsmath,amssymb,amsfonts,amsthm,graphicx}

\newtheorem{theorem}{Theorem}[section]

\newtheorem{corollary}[theorem]{Corollary}
\newtheorem{conjecture}[theorem]{Conjecture}
\numberwithin{equation}{section}

\theoremstyle{remark}

\newcommand{\Ric}{\mathop{\mathrm{Ric}}}

\def\Lap{\triangle}

\title{A Harnack inequality for the parabolic Allen-Cahn equation}

\author{Mihai B\u{a}ile\c{s}teanu}
\thanks{Department of Mathematical Sciences, Central Connecticut State University, 120 Marcus White Hall, New Britain, CT 06050, USA \texttt{mihaib@ccsu.edu}}

\begin{document}

\begin{abstract}
We prove a differential Harnack inequality for the solution of the parabolic Allen-Cahn equation $ \frac{\partial f}{\partial t}=\Lap f-(f^3-f)$ on a closed n-dimensional manifold. As a corollary we find a classical Harnack inequality. We also formally compare the standing wave solution to a gradient estimate of Modica from the 1980s for the elliptic equation.   
\end{abstract}

\maketitle

\section{Introduction}

The elliptic Allen-Cahn equation $\Lap f=f^3-f$, with $|f|\leq 1$ is a very popular non-linear PDE, which gave rise to an interesting conjecture formulated by De Giorgi in \cite{DeGiorgi79}:
\begin{conjecture}
Suppose that $f$ is an entire solution to the Allen-Cahn equation
\begin{align}\label{A-C-elliptic}
\Lap f-f^3+f=0, \hspace{1cm} |f|\leq 1, \hspace{1cm} x=(x',x_n)\in\mathbb{R}^n 
\end{align}
satisfying $\frac{\partial f}{\partial x_n}>0$ for $x\in\mathbb{R}^n$. Then, at least for $n\leq 8$, the level sets of $f$ must be hyperplanes. 
\end{conjecture}
Equivalently this can be stated as follows: any entire solution which is monotone in one direction should be one-dimensional. The conjecture was proved for $n=2$ by Ghoussoub-Gui \cite{ghoussoub-gui98}, for $n=3$ by Ambrosio-Cabre \cite{ambrosio-cabre00} and for $4\leq n\leq 8$ by Savin, under an extra assumption (see \cite{savin03,savin09,savin10}). The equation has been deeply investigated, see for example: \cite{gui97, chen92a, chen92b, chen94,chen04,chen07,chen06,chen00, kowalczyk12b,gui10,ghoussoub-gui03,gui03,gui12b,gui12a,kowalczyk12a,kowalczyk15}, but still offers many interesting questions. Note that De Giorgi's conjecture is false for $n\geq 9$ (see, for example, \cite{kowalczyk11}).

Originally, the equation appeared in the study of the process of phase separation in iron alloys, including order-disorder transitions (see \cite{allen-cahn79}). Moreover, it is related to the study of minimal surfaces, so it is an interesting topic for geometry too.

In this paper we are focusing on the parabolic Allen-Cahn equation:
\begin{align}\label{A-C-parabolic}
\begin{cases}
 \frac{\partial f}{\partial t}=\Lap f-(f^3-f) \\
 f(x,0)=f_0(x)>0
\end{cases}
\end{align}
where $f:\mathbb{R}^n\times[0,\infty)\to \mathbb{R}$.

Some solutions have interfaces that travel in one direction with a constant speed $c$. Without loss of generality, we may assume that the travel direction is $x_n$. If the solution is of the form $f(x,t)=p(x_1,...,x_{n-1}, x_n-ct)$, for some function $p:\mathbb{R}^n\to\mathbb{R}$, then $p(x)$ is a traveling wave solution with speed $c\geq 0$. Substituting $f$ with $p$ in (\ref{A-C-parabolic}) we obtain the following equation:
\begin{align}\label{A-C-travel}
c\frac{\partial p}{\partial x_n}+\Lap p =p^3-p 
\end{align}
In the case when $c=0$, if we impose the condition $\lim\limits_{x_n\to\pm\infty} p(x)=\pm 1$, for all $x\in\mathbb{R}^n$, then the traveling wave (in this case, a standing wave) satisfies the elliptic Allen-Cahn equation (\ref{A-C-elliptic}). Therefore, by obtaining estimates on the solution of the parabolic equation, it is possible to relate them to the solution of the elliptic one, by using the standing wave as an intermediary. For a nice introduction about traveling waves of the Allen-Cahn equation see \cite{fukao-morita-ninomiya04,morita-ninomiya10}.

The focus of this paper is to use geometric methods to obtain a Harnack inequality for the solution of the parabolic Allen-Cahn equation. The techniques only work for positive solutions and as such cannot be directly applied to get a bound on the elliptic equation. We nevertheless formally compare the standing wave solution to Modica's gradient estimates and show that they have similar polynomial bounds.    

The main result of this paper is the following theorem, which is a differential Harnack inequality for the solution of equation (\ref{A-C-parabolic}) on a closed manifold:
\begin{theorem}\label{Harnack-theorem}
Let $M$ be a compact Riemannian manifold without boundary of dimension $n$. Assume the Ricci curvature $\Ric\geq -k$, where $k\geq 0$. Consider $f$ to be a solution to the parabolic Allen-Cahn equation (\ref{A-C-parabolic}), such that $0<f<1$. Denote with $u(x,t)=\log f(x,t)$. Define the Harnack quantity $h(x,t)$ to be
\[h(x,t)=\Lap u+\alpha |\nabla u|^2+\beta e^{2u}+\phi(t)\]
where $\alpha\in(0,1)$ and $\beta\leq\min\left\{-\frac{n(\alpha+2)}{2\alpha^2-2\alpha+3n}, -\frac{nk}{4\alpha}, -\frac{n}{2(1-\alpha)}\right\}$ and  $\phi(t)=\frac{-2q\coth(-2qt)+b}{2a}$ ($a,b,d$ are positive real numbers that only depend on $\alpha$, $\beta$). Then $h(x,t)\geq 0$ for any $t>0$. 
\end{theorem}

This is a Li-Yau-Hamilton type Harnack inequality, the study of which began in the seminal paper by Li-Yau \cite{ly86}, and which has played a strong role in the study of geometric flows (see, for example, \cite{H93harnack, perelman1}). In the present paper we apply a method previously used in the study of Ricci flow and other heat-type equations (see, for example, the work of Cao \cite{cao08, cao-hamilton09}). The method both determines the Harnack quantity and proves that it's non-negative, using the maximum principle. The method has proven to be able to recover classical results for other non-linear parabolic equation (see \cite{cao-cerenzia-kazaras14, cao-ljun-liu13}) or for the curve shortening flow (\cite{bailesteanu15}).  The procedure has its roots in the study of Ricci flow, but recently its efficiency has been observed for other heat-type equations also (see, for example, \cite{cao-ljun-liu13} for a non-linear heat equation, or \cite{cao-cerenzia-kazaras14} for the Endangered Species Equation).

By integrating the above along a space-time path, we obtain the following classical Harnack inequality:

\begin{theorem}\label{classical-Harnack-theorem}
Let $M$ be a compact Riemannian manifold without boundary of dimension $n$. Assume the Ricci curvature $\Ric\geq -k$, where $k\geq 0$. Consider $f$ to be a solution to the parabolic Allen-Cahn equation (\ref{A-C-parabolic}), such that $0<f<1$. Then for $\alpha,\beta$ satisfying the conditions of theorem \ref{Harnack-theorem}: 
\[\frac{f(x_2,t_2)}{f(x_1,t_1)}\geq e^{-\frac{d(x_1,x_2)^2}{4(1-\alpha)(t_2-t_1)}}\cdot\left(\frac{e^{-2qt_2}-1}{e^{-2qt_1}-1}\right)^{-1/a}\cdot e^{-\frac{4q+b}{2a}(t_2-t_1)}\]
for $a,b,q$ positive real numbers that only depend on $\alpha$, $\beta$. 
\end{theorem}

The standing wave of this equation cannot remain in the interval $(0,1)$. This can be seen by comparing it to a radial solution in a ball with expanding radius and zero boundary condition, which is below the standing wave solution. By maximum principle the radial solution would remain below the standing wave, but as the radius goes to infinity, the radial solution would tend to $1$, thus the standing wave would also be $1$. Therefore the Harnack inequality cannot be used directly to estimate the elliptic solution. But we can at least show that the geometric method produces a similar polynomial bound as the classical PDE method. We will formally assume that the standing wave solution is in $(0,1)$ and obtain the following:  

\begin{corollary}\label{corollary-travel-wave}
Let $p$ be a standing wave solution for the parabolic Allen-Cahn equation (\ref{A-C-parabolic}) in $\mathbb{R}^n$ such that $0<p<1$. Then $p$ satisfies:
\begin{align}\label{stand-wave}
|\nabla p|^2\leq p^2[(2n-1)-(n-1)p^2] 
\end{align}
\end{corollary}

This is reminiscent of Modica's gradient estimate (\cite{modica85}):
\begin{theorem}\label{Modica}
Let $F\in C^2(\mathbb{R})$ be a non-negative function and $u\in C^3(\mathbb{R}^n)$ be a bounded entire solution to the equation $\Lap u=f(u)$, where $f=F'$. Then $|\nabla u|^2(x)\leq 2 F(u(x))$ for every $x\in\mathbb{R}^n$. 
\end{theorem}
In our case, Modica's result would translate as $|\nabla p|^2\leq \frac{1}{2}p^2(p^2-2)+\frac{1}{2}$ (note that a constant is added in order to obtain a non-negative anti-derivative of $x^3-x$). Comparing the two gradient estimates, one can notice that both bounds have a degree 4 polynomial and while the present result depends on the dimension of the underlying space, it does offer a polynomial ($x^2(2n-1-(n-1)x^2)$) that does not need to be non-negative everywhere. Moreover, as it will be shown in the last section, the new bound is an improvement on parts of the interval of solutions. 

The paper is structured as follows: in section 2 we describe the setting of the problem and summarize the procedure we use to construct the Harnack quantity. In section 3 we determine the exact expression of the Harnack quantity and show that it's non-negative, thus proving theorem \ref{Harnack-theorem}, while in section 4 we prove the classical Harnack inequality \ref{classical-Harnack-theorem}. The corollary \ref{corollary-travel-wave} involving the traveling wave and the comparision with Modica's result is proven in section 5. 

\textbf{Acknowledgements} The author would like to thank Prof. Xiaodong Cao for suggesting applying this method to the Allen-Cahn equation and for fruitful discussions.

\section{Background and setting}

We consider an $n$-dimensional orientable complete compact Riemannian manifold $M$, without boundary and with Ricci curvature $\Ric\geq -k$ for some $k\geq 0$. Let $f:M\times[0,\infty)\to \mathbb{R}$ be a solution to the parabolic Allen-Cahn equation: 
\begin{align}\label{A-C-parabolic-bis}
\begin{cases}
 \frac{\partial f}{\partial t}=\Lap f-(f^3-f) \\
 f(x,0)=f_0(x)>0
\end{cases}
\end{align}

Moreover, we assume that that $f_0(x)<1$. By maximum principle, if the function $f\in(0,1)$ at the beginning, then it will stay in $(0,1)$ at any time $t>0$ (this is because $f-f^3=f(1-f^2)$ stays non-negative for $0<f<1$). Therefore, $f(x,t)\in(0,1)$ from now on, for any $x\in M$ and $t\geq 0$. 

The method that we will use to determine the Harnack quantity and prove that it's non-negative is summarized as follows: let $f>0$ be the function for which one needs to prove a Harnack inequality, where $f$ satisfies a heat-type equation. If  $u=\log f$, then $u_t=\triangle u+|\nabla u|^2$+other terms. Inspired by the expression of $u_t$, one can define the Harnack quantity $h=\alpha\triangle u+\beta|\nabla u|^2+\text{possibly other terms}+\phi(x,t)$. $\phi$ is a function that has to go to infinity as $t\to 0$ and will be determined at the end. Assuming there is a first point (first, with respect to time) when $h(x,t)\leq 0$.  As $\phi$ is very large at time $0$, it means that  $\lim\limits_{t\to 0^+} h=\infty$, i.e. $h$ is positive close to $t=0$, so $t_1>0$. Moreover at this point, the time derivative of $h$ has to be negative. At that particular time, $h(x,t_1)$ is a local space minimum on the whole manifold, so if we apply the maximum principle, $\Lap h\geq 0$ and $\nabla h=0$. Therefore the quantity $h_t-\triangle h -2\nabla u\nabla h$ has to be non-negative (the fact that $M$ is closed allows for this to be done globally). We impose conditions on $\alpha, \beta$ and possibly other constants and restrictions for $\phi$ (from which one builds $\phi$) such that the expression turns out positive, contradicting the initial assumption. Therefore, we simultaneously determine $h$ and show that it's non-negative. 

We note here that it's possible to obtain a local Harnack estimate, in the case when $M$ is complete, non-compact, but this is more technical and it is being done in a subsequent paper.

\section{The Harnack quantity}

In this section we will determine the exact expression of the Harnack quantity and prove that it's positive for any time $t>0$. 

Let $f$ be a solution to (\ref{A-C-parabolic-bis}). Assuming that $0<f(x,y)<1$, we are able to introduce $u=\log f$ (and hence $-\infty<u<0$). It follows that $u$ satisfies:
\begin{align}\label{u-Allen-Cahn}
u_t=\Lap u+|\nabla u|^2+1-e^{2u} 
\end{align}

Let $h(x,t)$ be the following generic Harnack quantity:
\begin{align*}
h(x,t):=\Lap u +\alpha|\nabla u|^2+\beta e^{2u}+\varphi(x,t)
\end{align*}
where $\varphi(x,t)=\phi(t)+\psi(x)$ ($\phi$ is a time-dependent function, while $\psi$ is space-dependent). In the present situation, since we are working on a compact manifold, we won't need the spatial function. However, we keep it for the time being, for completeness of the method. 

Our goal is to compute the time evolution of $h$ and then use the maximum principle to establish the positivity of $h$ given a particular choice of $\alpha$, $\beta$ and $\phi(t)$. Moreover, $\phi(t)$ has to be very large at $t=0$, in order to dominate the other terms and ensure that at time close to $0$ the Harnack quantity is positive. 

Applying the heat operator to each term of $h(x,t)$, one obtains:
\begin{align*}
(\partial_t-\Lap)u &= |\nabla u|^2+1-e^{2u} \\
(\partial_t-\Lap)\varphi & = \phi_t-\Lap\psi \\
(\partial_t-\Lap)(\beta e^{2u}) & = 2\beta e^{2u}-2\beta e^{4u}-2\beta e^{2u}|\nabla u|^2=2\beta(e^{2u}-e^{4u}-e^{2u}|\nabla u|^2) \\
(\partial_t-\Lap)(\Lap u) & = \Lap |\nabla u|^2-2(\Lap u)e^{2u}-4|\nabla u|^2e^{2u}\\
(\partial_t-\Lap)(\alpha|\nabla u|^2) &= 2\alpha\nabla u\cdot \nabla(\Lap u)+2\alpha\nabla u\cdot\nabla(|\nabla u|^2)-4\alpha|\nabla u|^2e^{2u}-\alpha\Lap(|\nabla u|^2)\\
2\nabla u\cdot\nabla h & = 2\nabla u\cdot\nabla(\Lap u)+2\alpha\nabla u\cdot \nabla|\nabla u|^2+4\beta|\nabla u|^2e^{2u}+2\nabla u\cdot\nabla\psi
\end{align*}

Recall Bochner's identity:
\[\Lap|\nabla u|^2=2|\nabla\nabla u|^2+2\nabla u\cdot\nabla(\Lap u)+2\Ric(\nabla u,\nabla u)\]

Putting everything together yields the following time evolution for $h(x,t)$:
\begin{align*}
(\partial_t-\Lap)h-2\nabla u\cdot\nabla h & = 2(1-\alpha)|\nabla\nabla u|^2-2(\Lap u)e^{2u}-|\nabla u|^2(4e^{2u}+4\alpha e^{2u}+6\beta e^{2u}) \\
  & +2(1-\alpha)\Ric(\nabla u,\nabla u) +2\beta(e^{2u}-e^{4u})\\ & +\phi_t-\Lap\psi-2\nabla u\nabla\psi
\end{align*}

Using the fact that $|\nabla\nabla u|^2\geq\frac{1}{n}(\Lap u)^2$ and $\Ric\geq -k$, one gets the following inequality:

\begin{align*}
(\partial_t-\Lap)h-2\nabla u\cdot\nabla h & \geq \frac{2}{n}(1-\alpha)(\Lap u)^2-2(\Lap u)e^{2u}-|\nabla u|^2(4e^{2u}+4\alpha e^{2u}+6\beta e^{2u}) \\
  & -2(1-\alpha)k|\nabla u|^2 +2\beta(e^{2u}-e^{4u})+\phi_t-\Lap\psi-2\nabla u\nabla\psi
\end{align*}

$\Lap u$ can be replaced with $h-\alpha|\nabla u|^2-\beta e^{2u}-\phi-\psi$, which leads to the expression:

\begin{align}\label{Harnack-time-evolution}
(\partial_t-\Lap)h-2\nabla u\cdot\nabla h \geq &   \nonumber\\
 & h\left[\frac{2(1-\alpha)}{n}h -\frac{4(1-\alpha)}{n}(\alpha|\nabla u|^2+\beta e^{2u}+\phi+\psi)-2e^{2u}\right] \nonumber\\
 & +\left[\frac{2(1-\alpha)}{n}(\alpha^2|\nabla u|^4+2\phi\psi)-2k(1-\alpha)|\nabla u|^2+\frac{4\alpha(1-\alpha)}{n}\phi|\nabla u|^2\right. \nonumber \\
 & +\left.|\nabla u|^2e^{2u}\left(\frac{4\alpha\beta(1-\alpha)}{n}-6\beta-2\alpha-4\right)\right] \\
 & \left[e^{4u}\cdot\frac{2\beta^2(1-\alpha)}{n}+e^{2u}\left(\frac{4\beta(1-\alpha)}{n}\phi+2\phi+2\beta\right)+\right. \nonumber\\
 & \left.+\frac{2(1-\alpha)}{n}\phi^2+\phi_t\right] \nonumber\\
 & +\left[\frac{4\alpha(1-\alpha)}{n}\cdot \psi|\nabla u|^2-2\nabla u\nabla\psi+e^{2u}\psi\left(\frac{4\beta(1-\alpha)}{n}+2\right)\right.\nonumber \\
 &\left. +\frac{2(1-\alpha)}{n}\psi^2-\Lap\psi \right] \nonumber
\end{align}

Heuristically, we group terms that involve $|\nabla u|^2$, powers of $e^{2u}$ and $\psi$ separately. To make the computation easier to follow, we will denote each expression as follows:
\[P_1=\frac{2(1-\alpha)}{n}h -\frac{4(1-\alpha)}{n}(\alpha|\nabla u|^2+\beta e^{2u}+\phi+\psi)-2e^{2u}\]

\begin{align*}
P_2& =\frac{2(1-\alpha)}{n}(\alpha^2|\nabla u|^4+2\phi\psi)-2k(1-\alpha)|\nabla u|^2+\frac{4\alpha(1-\alpha)}{n}\phi|\nabla u|^2\\
 & +|\nabla u|^2e^{2u}\left(\frac{4\alpha\beta(1-\alpha)}{n}-6\beta-2\alpha-4\right)
\end{align*}

% (notice that the terms involving $4\alpha e^{2u}$ have canceled)

\[P_3=e^{4u}\cdot\frac{2\beta^2(1-\alpha)}{n}+e^{2u}\left(\frac{4\beta(1-\alpha)}{n}\phi+2\phi+2\beta\right)+\frac{2(1-\alpha)}{n}\phi^2+\phi_t\]

\[P_4=\frac{4\alpha(1-\alpha)}{n}\cdot \psi|\nabla u|^2-2\nabla u\nabla\psi+e^{2u}\psi\left(\frac{4\beta(1-\alpha)}{n}+2\right)+\frac{2(1-\alpha)}{n}\psi^2-\Lap\psi\]

We have thus shown that 
\[(\partial_t-\Lap)h-2\nabla u\cdot\nabla h \geq h P_1+P_2+P_3+P_4.\]

Assume there is a first time when $h\leq 0$. Since $\varphi(x,t)$ will be constructed such that it goes to infinity at time $0$, $h$ has to be positive close to the starting time. But since the solution is smooth (it's a heat-type equation), $h$ is also smooth, so there has to be a first time $t_0>0$ when $\min_M h(x,t)=0$. Recall that $M$ is a compact manifold, so by applying the maximum principle $h_t\leq 0$, $\nabla h=0$ and $\Lap h\geq 0$ at that time (the function $h(x,t_0)$ has a spatial minimum on the whole $M$ at $(x_0,t_0)$). Moreover, at this point $\Lap u=-\alpha |\nabla u|^2-\beta e^{2u}-\phi-\psi$. Therefore at $(x_0,t_0)$:
\[0\geq P_2+P_3+P_4.\]

We want to obtain a contradiction, so all we need to do is find specific values for $\alpha,\beta$ and $\varphi$ such that $P_2+P_3+P_4>0$. 

In the present situation, the manifold is closed, so we don't have to worry about boundary issues (the minimum will be spatially global), therefore we don't need the spatial component for $\varphi$, hence $\psi=0$ or $\varphi(x,t)=\phi(t)$.  Thus, $P_4=0$. 

All we are left to do is find $\alpha,\beta, \phi(t)$ such that $P_2$ and $P_3$ are non-negative and at least one of them positive. 

We start by analyzing $P_2$. Recall that 
\begin{align*}
P_2& =\frac{2(1-\alpha)}{n}\alpha^2|\nabla u|^4-2k(1-\alpha)|\nabla u|^2+\frac{4\alpha(1-\alpha)}{n}\phi|\nabla u|^2\\
 & +|\nabla u|^2e^{2u}\left(\frac{4\alpha\beta(1-\alpha)}{n}-6\beta-2\alpha-4\right)
\end{align*}

Notice that since $\psi=0$, the term involving $\phi\psi$ disappeared. Imposing that $\alpha\in(0,1)$ assures that the first is non-negative. To ensure the positivity of the sum of the second and third term, we need $\phi\geq \frac{nk}{2\alpha}$. Finally, for the last term to be non-negative, we need $\frac{4\alpha\beta(1-\alpha)}{n}-6\beta-2\alpha-4\geq 0$, which is equivalent to 
\[\beta\leq -\frac{n(\alpha+2)}{2\alpha^2-2\alpha+3n}<0\]

Note that this means $\beta$ has to be negative.

Next we analyze $P_3$. Recall that: 
\[P_3= e^{4u}\cdot\frac{2\beta^2(1-\alpha)}{n}+2e^{2u}\left[\left(\frac{2\beta(1-\alpha)}{n}+1\right)\phi+\beta\right]+\frac{2(1-\alpha)}{n}\phi^2+\phi_t\]

By applying the inequality $p x^2+2q x\geq -\frac{q^2}{p}$ for $x=e^{2u}$, we obtain:

\begin{align*}
P_3&\geq -\frac{\left[\left(\frac{2\beta(1-\alpha)}{n}+1\right)\phi+\beta\right]^2}{\frac{2\beta^2(1-\alpha)}{n}}+\frac{2(1-\alpha)}{n}\phi^2+\phi_t\\
   &=-\frac{2}{\beta}\left(1+\frac{n}{4\beta(1-\alpha)}\right)\phi^2-2\left(1+\frac{n}{2\beta(1-\alpha)}\right)\phi-\frac{2n}{1-\alpha}+\phi_t
\end{align*}

At this point, we need to solve a differential inequality, and there are two situations: first, if $\beta\geq -\frac{n}{4(1-\alpha)}$ (which is smaller than $-\frac{n(\alpha+2)}{2\alpha^2-2\alpha+3n}$ for any $n\geq 3$), then both expressions $\left(1+\frac{n}{2\beta(1-\alpha)}\right)$ and $\left(1+\frac{n}{4\beta(1-\alpha)}\right)$ are negative. This means that we want $\phi$ to satisfy:
\[-a \phi^2+b\phi -c+\phi_t\geq 0\] for some positive numbers $a,b,c>0$.

Second, if $\beta\leq -\frac{n}{2(1-\alpha)}$ (which is also smaller than $-\frac{n(\alpha+2)}{2\alpha^2-2\alpha+3n}$ for any $n\geq 1$ and $\alpha\geq \frac{1}{2}$ or for any $n\geq 2$) then both of those expressions are positive, so $\phi$ would satisfy:
\[a \phi^2-b\phi -c+\phi_t\geq 0\] for some positive numbers $a,b,c>0$.

Since the second case gives a solution for a larger range of values for $\beta$ and $n$, we will focus on it. For an arbitrarily small $d>0$, if we determine $\phi>0$ such that 
\[a \phi^2-(b+d)\phi -c+\phi_t=0\] then we will ensure that $P_3>0$.

The above equation has an exact solution, given by: 
 \[\phi_d(t)=\frac{-2q\coth(-qt)+(b+d)}{2a}=\frac{[(b+d)-2q]-e^{2qt}(b+d+2q)}{2a(1-e^{2qt})}\]
where $a=-\frac{2}{\beta}\left(1+\frac{n}{4\beta(1-\alpha)}\right)$, $b=2\left(1+\frac{n}{2\beta(1-\alpha)}\right)$, $c=\frac{2n}{1-\alpha}$ and $q=\sqrt{ac+\frac{1}{4}(b+d)^2}$. Notice that, in fact, $c$ depends on $a$ and $b$ since $2+\frac{c}{2\beta}=b$ and $2+\frac{c}{4\beta}=-a\beta$. Moreover, $q$ is well defined, since $a,c>0$ and thus the expression under the root is positive. 

The expression of $\phi_d(t)$ from above is positive for all $t$ and its limit, as $t\to 0^+$ is $\infty$ ($a,b,q,c,d>0$). These are the desired properties for $\phi_d(t)$ that guarantee that $P_3>0$. 
 
Recall that in order for $P_2>0$, $\phi(t)$ must satisfy $\phi\geq \frac{nk}{2\alpha}$, which in turns imposes the condition $\frac{b+d+2q}{2a}\geq \frac{nk}{2\alpha}$. By replacing $q$ with the expression above, we get the following condition:
\[\frac{b+d+2\sqrt{ac+\frac{1}{4}(b+d)^2}}{2a}\geq\frac{nk}{2\alpha}.\]

Since $ac$ is positive:  
\[\frac{b+d+2\sqrt{ac+\frac{1}{4}(b+d)^2}}{2a}\geq\frac{b+d}{2a}.\]

If one imposes that $\frac{b}{a}\geq \frac{nk}{2\alpha}$, since $d,a$ are positive, then the condition on $\phi$ is satisfied. We will show that this leads to an upper bound on $\beta$. 

First, $\frac{b}{a}=-\beta\frac{2\beta(1-\alpha)+n}{4\beta(1-\alpha)+n}\geq -2\beta$ since the quotient is less than $1$, which follows from the fact that $\beta<0$. Therefore, we only need the condition that $-2\beta\geq \frac{nk}{2\alpha}$, which is equivalent to saying that $\beta\leq -\frac{nk}{4\alpha}$. 
 
In conclusion, for $\alpha\in(0,1)$, $\beta\leq\min\{-\frac{n(\alpha+2)}{2\alpha^2-2\alpha+3n}, -\frac{nk}{4\alpha}, -\frac{n}{2(1-\alpha)}\}$ and  $\phi_d(t)$ having the above expression (with $a,b,d,q$ only depending on $\alpha$ and  $\beta$), $P_2+P_3>0$. This contradicts the fact that at the point $(x_0,t_0)$, $P_2+P_3\leq 0$, hence there does not exist any point for which $h\leq 0$, hence $h_d(x,t)= \Lap u +\alpha|\nabla u|^2+\beta e^{2u}+\phi_d(x,t)> 0$ for any $t>0$.

The final step to prove theorem \ref{Harnack-theorem} is to take $d\to 0$, which will turn $\phi_d(t)$ into $\phi(t)=\frac{-2q\coth(-qt)+b}{2a}=\frac{b-q-e^{2qt}(b+2q)}{2a(1-e^{2qt})}$, where $q=\frac{1}{2}\sqrt{b^2+4ac}$.

\section{The Classical Harnack inequality} 

In this section we will obtain a classical Harnack inequality. The procedure will follow the classical method developed by Li-Yau and consists of integrating the differential Harnack estimate along a space time-curve.  

We choose two points $(x_1,t_1)$ and $(x_2,t_2)$ and a space-time path connecting them $\gamma:[t_1,t_2]\to M$, such that $\gamma(t_1)=x_1$ and $\gamma(t_2)=x_2$. The value of the function $u(x,t)=\log f(x,t)$ along the path is given by $v(t):=u(\gamma(t),t)$. Its time derivative has the expression $v'(t)=u_t+\nabla u\cdot\frac{d\gamma}{dt}$, which is the same as
\[v'(t)=\Lap u+|\nabla u|^2+1-e^{2u}+\nabla u\cdot\frac{d\gamma}{dt}.\]
From the Harnack inequality, one can bound $\Lap u$ and get that $\Lap u\geq -\alpha|\nabla u|^2-\beta e^{2u}-\phi(t)$ which gives a lower bound for $v'(t)$:
\[v'(t)\geq (1-\alpha) |\nabla u|^2-(\beta+1) e^{2u}+1-\phi(t)+\nabla u\cdot\frac{d\gamma}{dt}.\]
Recall that $\alpha\in(0,1)$ and $\beta\leq\min\left\{-\frac{n(\alpha+2)}{2\alpha^2-2\alpha+3n}, -\frac{nk}{4\alpha}, -\frac{n}{2(1-\alpha)}\right\}$.

Using the inequality $ay^2+by\geq-\frac{b^2}{4a}$ for $y=|\nabla u|$, we obtain that 
\[v'(t)\geq -\frac{1}{4(1-\alpha)}\left|\frac{d\gamma}{dt}\right|^2-(\beta+1)e^{2u}+1-\phi(t).\]

Since $f(x,t)\in(0,1)$, it follows that $e^{2u}=f^2\in(0,1)$, so $-e^{2u}\geq -1$. Moreover, $\beta$ is negative, so $-\beta e^{2u}\geq 0$, therefore one finally has:
\[v'(t)\geq -\frac{1}{4(1-\alpha)}\left|\frac{d\gamma}{dt}\right|^2-\phi(t).\]

First, we integrate in time from $t_1$ and $t_2$ and we obtain
\[v(t_2)-v(t_1)=\int\limits_{t_1}^{t_2}v'(t)\ dt\geq -\frac{1}{4(1-\alpha)}\int\limits_{t_1}^{t_2}\left|\frac{d\gamma}{dt}\right|^2\ dt-\int\limits_{t_1}^{t_2}\phi(t)\ dt.\]

If one chooses $\gamma$ to be a minimizing geodesic between the endpoints, the first integral becomes $\frac{d(x_1,x_2)^2}{t_2-t_1}$. 

The second integral can be computed directly:
\[-\int\limits_{t_1}^{t_2}\phi(t)\ dt=-\int\limits_{t_1}^{t_2}\frac{-2q\coth(-qt)+b}{2a}\ dt=-\frac{4q+b}{2a}(t_2-t_1)-\frac{1}{a}\ln\left(\frac{e^{-2qt_2}-1}{e^{-2qt_1}-1}\right)\]

with $a=-\frac{2}{\beta}\left(1+\frac{n}{4\beta(1-\alpha)}\right)$, $b=2\left(1+\frac{n}{2\beta(1-\alpha)}\right)$, $c=\frac{2n}{1-\alpha}$ and $q=\frac{1}{2}\sqrt{b^2+4ac}$.

Finally, by noticing that $v(t_2)-v(t_1)=u(x_2,t_2)-u(x_1,t_1)=\ln\frac{f(x_2,t_2)}{f(x_1,t_1)}$ and taking the exponential we get the result of theorem \ref{classical-Harnack-theorem}.

If $M$ has non-negative Ricci curvature, one may choose $\alpha=\frac{1}{2}$ and $\beta=-n$. It follows that $a=1/n$, $b=0$, $c=4n$ and $q=2$. This turns the previous classical Harnack inequality into a simpler form:
\[\frac{f(x_2,t_2)}{f(x_1,t_1)}\geq e^{-\frac{d(x_1,x_2)^2+8n(t_2-t_1)^2}{2(t_2-t_1)}}\left(\frac{e^{-4t_2}-1}{e^{-4t_1}-1}\right)^{-n}\]

\section{The traveling wave solution}

We will now focus on the traveling wave solution, when when $M=\mathbb{R}^n$ (hence the Ricci curvature is $0$). Like at the end of the previous section, we choose $\alpha=\frac{1}{2}$ and $\beta=-n$, making $a=1/n$, $b=0$, $c=4n$ and $q=2$. The expression of $\phi(t)$ becomes:
\[\phi(t)=2n\frac{e^{4t}+1}{e^{4t}-1}\]

Notice that $\lim\limits_{t\to\infty}\phi(t)=2n$.

Theorem \ref{Harnack-theorem} says that the quantity $h(x,t)=\Lap u+\frac{1}{2}|\nabla u|^2-n e^{2u}+\phi(t)$ is always non-negative. Replacing $u$ with $\log f$, the above becomes:
\[\frac{f_t}{f}-\frac{1}{2}\frac{|\nabla f|^2}{f^2}+(1-n)f^2+\phi(t)-1\geq 0.\]

Assume that $f$ is in fact a traveling wave with speed $c\geq 0$, i.e. there is a function $p:\mathbb{R}^n\to\mathbb{R}$ such that $f(x,t)=p(x_1,x_2,...,x_n-ct)$. In that case the above inequality becomes:
\[-c\frac{\partial_n p}{p}-\frac{|\nabla p|^2}{2p^2}+(1-n)p^2+\phi(t)-1\geq 0\]

Finally, in the case when the traveling wave is in fact stationary, i.e. $c=0$ (recall that, in that case, $p$ satisfies the elliptic Allen-Cahn equation (\ref{A-C-elliptic})), the above inequality becomes:
\[-\frac{|\nabla p|^2}{p^2}+(1-n)p^2+\phi(t)-1\geq 0\]
Since the wave is stationary, it will look the same at any time, so the inequality is true as $t\to\infty$, which leads to:
\[-\frac{|\nabla p|^2}{p^2}-p^2(n-1)+2n-1\geq 0\]
Therefore we obtain
\[|\nabla p|^2\leq p^2[2n-1-(n-1)p^2]\]  
proving Corollary \ref{corollary-travel-wave}. 

\bigskip
\textbf{Comparison with Modica's result} Recall that Modica's result asserts that a solution of the equation $\Lap u=F'(u)$ satisfies $|\nabla u|^2\leq 2 F(u)$, for a non-negative function $F$. In the Allen-Cahn equation case, this translates to $|\nabla u|^2\leq \frac{1}{2}u^2(u^2-2)+\frac{1}{2}$. A visual comparison between the two degree 4 polynomials suggests that our estimate would be an improvement around $0$, if proven to hold for $f\in[-1,1]$. This is because on the interval $[-1,1]$ the graph of $x^2(3-x^2)$ sits mostly below the graph of $\frac{1}{2}x^2(x^2-2)+\frac{1}{2}$ (we consider $n=2$, and notice that for larger $n$ the estimate is worse).

\begin{center}
\includegraphics[scale=0.3]{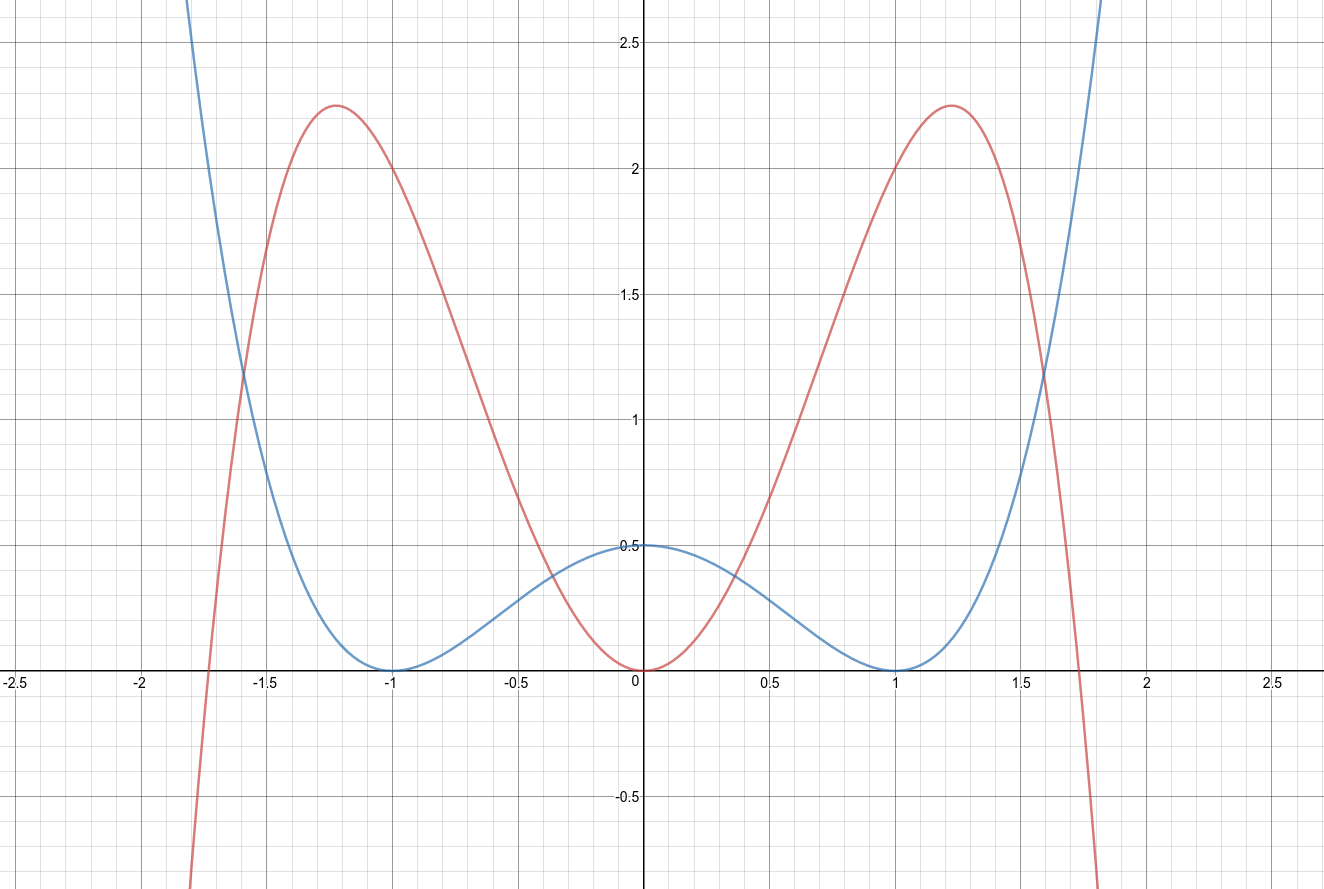} 

Figure 1: Comparision between $\frac{1}{2}x^2(x^2-2)+\frac{1}{2}$ (blue) and $x^2(3-x^2)$ (red)
\end{center}

\bibliographystyle{abbrv}
\bibliography{bio}

\end{document}